\documentclass[11pt,a4paper]{article}

\usepackage{amsmath,amssymb,amsfonts,amsthm,latexsym,tikz}
\usepackage[margin=1in]{geometry}
\usepackage{graphicx, bm, bbm, yfonts, colortbl}
\usepackage[citecolor=blue,urlcolor=blue,hyperindex,breaklinks]{hyperref}
\usepackage{subcaption,float}
\usepackage{MnSymbol}
\usepackage{verbatim}

\newcommand{\beq}{\begin{equation}}
\newcommand{\eeq}{\end{equation}}
\newcommand{\beqn}{\begin{equation*}}
\newcommand{\eeqn}{\end{equation*}}

\numberwithin{equation}{section}
\numberwithin{subsection}{section}

\begin{document}
\begin{center}	
	{ \scshape \bf \large Does 2026 AI exhibit intelligence, \\
	 or can Claude outsmart Pierre or Catherine ?}\footnote[1]{Preliminary version, May 26, 2026.
	 

}

	\vspace{0.5 in} 
	
	{\scshape {\large Robert C.~Dalang}\\ \'Ecole Polytechnique F\'ed\'erale de Lausanne} 
	\vspace{0.5\baselineskip}

	\vspace{0.5 in} 

\end{center}

\begin{abstract}
Using a sequence of high-school level mathematics questions that were not available on the Internet, we compare the performance of the popular AI software Claude with that of my friends and fellow human beings Pierre and Catherine. Pierre had solid scientific training as a young man, while Catherine studied literature. All three were subjected to a simulated pre-calculus oral exam with main questions and follow-up questions. Their performances are compared and the ones with the best and worst performances are identified. The outcome is that the current version of Claude, even though it is an extremely useful tool that has probably recorded the solution to nearly all calculus questions that are available on the Internet, {\em exhibits only a very limited understanding of the subject} and {\em does not exhibit the ability to make intelligent connections} between different features of a pre-calculus mathematics problem that it has never seen before.

\end{abstract}

\section{Introduction}\label{s1}

Just asking the question posed in the title of this article might suggest that the author knows what intelligence is. Of course, as for most of the important questions in life (What is democracy? What is free choice ? What am I looking for in a future life partner?), we do not really know what we are asking about! But as my mother used to say, ``you know it when you see it,'' and maybe the same old-fashioned common sense can be applied to Artificial Intelligence (AI) in the year 2026 and to intelligence. Another option is to look up the meaning of {\em intelligence} in your favorite dictionary or online ressource. My old Webster dictionary begins with ``Capacity for reasoning, understanding, and for similar forms of mental activity; aptitude in grasping truths, facts, meanings, etc.'' (it also mentions {\em military intelligence,} which is not under discussion here). Ressources on the web give lengthier descriptions. I consider that we really want to distinguish between knowledge and intelligence: intelligence is more ``putting two and two together'' than ``knowing that $2+2=4$,'' it is also more ``connecting the dots'' than ``knowing where the dots are.''
 
I have forty years of experience teaching mathematics at the university level, and I have seen my fair share of bright students who brilliantly pass all their exams, as well as many dedicated students who just cannot make it and unfortunately leave without their degree. I do not claim that anything that I have taught or the exams that I have given directly tests a student's ``intelligence,'' but mathematics is often used by some institutions of higher learning to select ``the best'' students, so there may be a relationship between some forms of intelligence and mathematical ability. And given that many humans beings consider that they are more intelligent than other animals, and that we have never seen other animals memorize the digits of $\pi$, this mathematical ability  could be some form of confirmation of the quality that we are considering.
 
The typical way of deciding if a student has passed the course and can proceed to the next one is by subjecting the student to an exam. There are many different kinds of exams: written exams with open questions, multiple choice exams and oral exams. At most European universities, oral exams used to be the preferred method for testing students, even for large classes (this method was done away with when the numbers of students increased substantially, faculty members preferred using their time to apply for grants and more time was needed for quality control of teaching). But those of us who have extensively used oral exams generally agree that, despite various weaknesses, it is a very good way of determining the student's understanding of the course materials and whether or nor he or she has mastered the required skills. Indeed, when a student has difficulty answering one question, the examiner follows up with a second question, and sometimes even with a sequence of questions, during which the examiner can probe the student's understanding of the material and determine fairly precisely the student's level of understanding and skill.

There are many different kinds of exam questions, and much thought has been given to what constitutes a ``good'' exam question. For instance, if the objective is to guarantee that every year, the exam has exactly the same level of difficulty, a first idea would simply be to give the same exam questions every year. Of course, this method fails completely, because students have a good memory, and this year's students will record the questions for next year's students, who can memorize the answers just before the exam, and then the course and the exam no longer serve any purpose. We will see that there is a similar problem when attempting to test AI, because AI does indeed have a very good and very long memory!

A different idea is that the exam questions should be different every year, and {\em close, but not identical,} to material and questions that have been discussed during the course. Indeed, we do not want students to learn questions and answers by heart (except in the rare case where the purpose is to test memorization abilities), and then simply reproduce the answers that were memorized before the exam.

Ideally, you could imagine that there is some very large set of all possible questions, and some metric on this set, that is, a method for determining ``how different'' questions are from one another, and then the exam should consist of questions that are not too far away from questions that have been discussed in the course, solvable with the knowledge that the students have acquired during the course, but not identical to the questions explicitly discussed in the course nor to the questions posed in previous years.

And you can well imagine that after decades of teaching calculus, by hundreds of professors around the world, each who posts his/her questions and solutions on the web, this strategy gets harder and harder to implement every year!

Fortunately for professors, students simply do not have the time to look at all the calculus exam questions that have ever been posed somewhere in the world at some point in time, and so, even if a question has been given before, somewhere and sometime, it is quite unlikely that a student will have seen it and very likely that the student will actually have to find the answer him/herself by using the knowledge acquired during the course. 

Unfortunately for attempts at testing the calculus abilities of current AI softwares, these Large Language programs have indeed had the time and ressources to read and digest every calculus question that has ever been given somewhere in the world. 

For example, a reasonably good exam on a calculus course for students could be designed according to the following principles: a data-base consisting of twenty-thousand questions, with their solution, is made publicly available. The final exam consists in a random selection of twenty of these questions. The AI-approach to passing this exam is simply to read and store the twenty-thousand question and answers, and thus {\em it will obtain a perfect score on the test,} even if it has no understanding of calculus. A human student could attempt to approximate this method, by studying as many of the twenty-thousand questions as possible. However, most students will prefer to attend the semester-long course, {\em learn the dozen or so key concepts} that make it possible, using their innate intelligence, to solve all twenty-thousand questions, {\em practice on 40 to 60 such questions,} and then they are not only ready to pass the test with flying colors, but they have actually learned some calculus!

It is with these thoughts in mind that I fell by chance on some questions that could be used to test the calculus abilities of 2026 AI software. In fact, the questions are even at the level of high-school mathematics, at least in most European and North American countries (these are the only places where I have working experience), so they can be attempted and appreciated even by people who have little inclination for mathematics.

\section{The circumstances that led to the comparison of Claude, Pierre and Catherine}\label{sec2}
A couple of weeks ago, I wanted to use part of the time in my last class of the Spring 2026 semester to discuss with my students the many relationships between mathematics and reality. As part of this effort, I unearthed a pre-calculus-level question that I had designed about 30 years ago, while a professor at Tufts University in Boston, where I helped to prepare a practice exam for our first semester calculus students. One of the question sets that I came up with at the time was a bit unusual, probably in part because it was during a year when the university wanted to emphasize the importance of writing in English, and so faculty members in all subjects were encouraged to formulate homework and exam problems in such a way that students would actually write at least part of the answer in plain English, rather than just reply by using formulas and mathematical symbols. I recall that at the time, some faculty members commented on the fact that my questions were rather original, and nevertheless very suitable for the students. I imagine that once given, the practice exam was put in a drawer for a few years before being discarded. So the only place where this question set has existed since the birth of twenty-first century AI is in my own mind. I figured that this would give me {\em one shot} at testing it on an AI software, because once the question is made available, any AI software can easily be trained to solve just this type of question.

And another interesting fact is that given the nature of the question, I was able, during a single 24-hour period, on May 17 and 18, 2026, to pose the question both to a human being (and a few days later, to a second one), presumed intelligent, and to the AI software Claude (by Anthropic, that I have switched to recently, after using other popular AI software for about a year and a half). 

Indeed, I first posed the main questions to a friend, that I will call Pierre. Pierre is a geologist with a Ph.D., therefore scientifically very well trained, but the last time that he took a math course was approximately forty years ago. But he is of course still very interested in science and remembers lots of what he had learned. When I asked him the main questions, and during follow-up questions, you could ``see the cogs wheels in his brain churning'' so as to connect the bits of information that he was gathering from the original question and my follow-up questions. To me, this looked like ``intelligence at work.'' 

The next day, I posed the same main questions to Claude. The follow-up questions were of course different, because they are the examiner's reaction to answers that are incorrect, and there is no way of knowing in advance what responses the examiner will receive and what questions should be used for follow-up. This is exactly the same as in an oral exam, where two consecutive students, who do not have the chance to talk to each other since the first one exits the exam room just as the second one enters, may be given the same question, and then the examiner can rank the two students and each will receive the grade that the examiner feels they merit.

Finally, a few days later, I asked the same questions to Catherine. While at the university some 40 years ago, Catherine studied literature, and so the last time that she saw a calculus question was probably in high-school. So she obviously had less scientific training than Pierre, and far much less stored knowledge from the Internet than Claude.

\section{An unusual question at the level of high-school mathematics}\label{sec3}

\begin{figure}[h]
\begin{center}
\includegraphics[height=4cm]{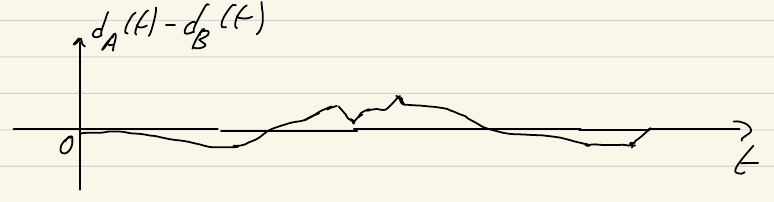}
\caption{The graph of a function. \label{fig1}}
 \end{center}
\end{figure}

The simulated oral exam and subsequent conversation takes place as follows. The responder (Claude,  Pierre or Catherine) is presented with the graph in Figure \ref{fig1}. Then the examiner explains that this graph shows the difference between the distance traveled by Swimmer $A$ and the distance travelled by Swimmer $B$, as a function of time, during an Olympic race. The responder can then examine the graph carefully.
Once the responder has had sufficient time to examine the graph, the examiner presents the following list of questions:
\medskip

   (a) What was the length of the race: (i) 50 meters, \hskip .3cm  (ii) 100 meters, \hskip .3cm or (iii) 200 m ?
\medskip

(b) Who is in the lead at the end of the first leg of the race ?
\medskip

(c) Who won the race ?
\medskip

(d) Imagine that you are a reporter for your favorite radio station, and are covering this race. How do you report the race to your listeners ?
\bigskip

I should clarify that the questions were asked in French and the conversation took place in French, which is the main language for Bachelor-level courses at my university. 
During the simulated exam, I treat Claude as though it were a human being taking the exam.

\section{Claude's exam}

Let's begin by recounting the exam with Claude. 
I must specify that the account presented here is an abridged and edited account, that aims to convey accurately the content of Claude's replies, as well as to show how the sequence of questions posed by the examiner is ajusted (in real time) {\em to zero-in on the actual understanding of the examinee} and his/her/its ability to piece together the relevant elements of acquired knowledge, versus simply evaluating the memorized knowledge that the examinee has acquired before the exam. Readers interested in the verbatim word-by-word ten-page transcript of my conversation with Claude can read \cite{claude}, which also describes how this transcript was produced.

Next, I should mention that since the questions are based on a picture, I was concerned that Claude would be at a disadvantage relative to Pierre and Catherine, who could actually ``see'' the picture, while Claude had to ``learn'' the picture by decoding the .pdf file that I gave it. However, during the conversation, it became clear that Claude could ``see'' the picture accurately, because it was able to point out explicitly all the relevant features of the graph, such as the local minima and maxima, the times where the curve crosses the axes, and the time where the curve finally returns to $0$ and stays there. Clearly, Claude was able to ``see'' the figure, determine and describe all of its important features, and even redraw the figure in a .pdf file of Claude's making, while adding to the figure the descriptive features that were part of Claude's answer. So I am convinced that Claude was able to acquire essentially the same knowledge about the curve as would a human being who looks at it for a few minutes. The {\em big question of interest,} however, is ``Can Claude, Pierre or Catherine {\em put this acquired knowledge to good use} so as to produce the correct interpretations?'' thereby exhibiting some form of understanding of the picture and some form of intelligence.

After telling Claude that I had prepared questions for engineering students (this question was in fact presented to my students as part of a larger list of questions at various levels of mathematics), I presented it with Question (a) concerning the length of the race. 

Claude's first response was: {\color{blue} ``The curve starts at $0$, dips slightly negative ($B$ briefly ahead), rises with oscillations ($A$ in front), and returns to 0 at the end---both swimmers finish the race at the same time (or nearly). The graph shows repeated oscillations in the middle, then a final descent. 

The oscillations are the key. These small oscillations correspond to the {\em tumble turns}---instants when one swimmer catches up or loses ground on the other at each half-length. I can count approximately ten half-lengths, meaning about 4 to 5 times the length of the pool. 

Conclusion: This is most likely a 200-meter race (4 × 50 m), which is the typical distance in which this kind of differential profile---with several visible turns and progressive race dynamics---is characteristic.
The expected answer is therefore: a 200 m race.'' }

There is some good and some bad in this response, but rather than comment on it to Claude, I decided to proceed with  Question (b): who was in the lead at the end of the first leg of the race ? 

Claude responded {\color{blue} ``After the first crossing, $d_A(t) - d_B(t)$ is negative, so $d_A < d_B$: swimmer $B$ is in the lead.''}

 Since this was not the answer that I was expecting, I then asked ``Can you show me on the graph the moment when the leading swimmer turns for the first time ?'' 
 
 Claude's reply:  {\color{blue} ``Of course! The first turn corresponds to the instant when the curve attains its first local minimum (when $B$, who is ahead, touches the wall).} 
 
\begin{figure}[h]
\begin{center}
\includegraphics[height=4cm]{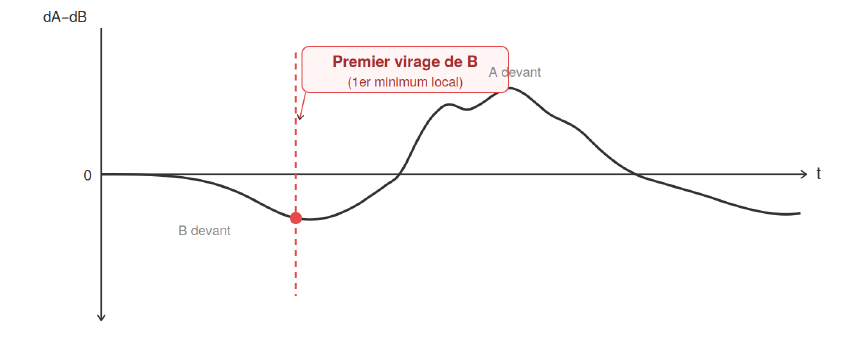}
\caption{The time at which the lead swimmer completes the first leg of the race, according to Claude. The speech bubble reads ``First turn of $B$ ($1^{\scriptsize \rm st}$ local minimum).'' Under the graph, it is stated ``B is ahead.'' \label{fig3}}
 \end{center}
\end{figure}

 And Claude drew on the graph a vertical  line that passed through the first local minimum and indicated that this was the time when the first swimmer turned around (see Figure \ref{fig3}). 
This answer is incorrect, and  the picture was even way off: for instance, his curve does not come back to the origin at the end. But I decided not to push this point just yet. I asked Claude to answer Question (d): how would you report the race to your radio listeners ? Claude then told the story of a race with several legs (but did not mention a precise number),  that ends in a {\em photo finish,} claiming (incorrectly) that the two swimmers are almost tied and that the winner will only be determined once the official times are given!

This answer did correspond somewhat to the fact that he thought that it was a 200 m race, but I was still disturbed about how Claude described the end of the race. So I asked Claude  to answer Question (c): who won the race ?

Claude replied {\color{blue} ``At the very end, $d_A(t) - d_B(t)$ is slightly negative: $d_A < d_B$, so $B$ won, by a very small margin.'' }

This answer is partly correct, since $B$ did win, but I was surprised that he said that $B$ won ``by a very small margin.'' So I decide to probe this issue and I asked it: ``Can you show me on the curve the moment when $B$ wins the race?'' 

\begin{figure}[h]
\begin{center}
\includegraphics[height=4cm]{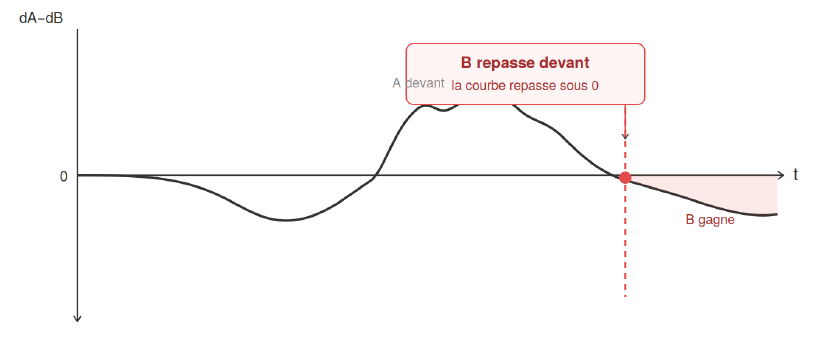}
\caption{The time at which Swimmer $B$ takes the lead and stays there, according to Claude. The speech bubble reads ``$B$ passes in front, the curve goes negative.'' \label{fig4}}
 \end{center}
\end{figure}

Claude replied {\color{blue} ``The key point is the second crossing of zero near the end of the race:''} and produced the picture in Figure \ref{fig4}.
It even added: {\color{blue} ``This is the time at which $B$ passes for the final time in front of $A$. The red zone shows the final meters where $d_A - d_B < 0$, i.e.~$B$ is in front through to the finish.''}

This reply is correct, but it made me realize that Claude was actually telling me at what point $B$ took the lead and kept it until the end of the race, rather than showing me the precise point in time where $B$ finished and finally won the race. My question could have been, by a small stretch,  interpreted in this way, so I decided to revisit the fact that his curve did not return to $0$ and stay there. I stated ``On my graph, the last point of the curve is at $0$: $d_A - d_B = 0$. That is not the case on your graph.'' 

\begin{figure}[h]
\begin{center}
\includegraphics[height=4cm]{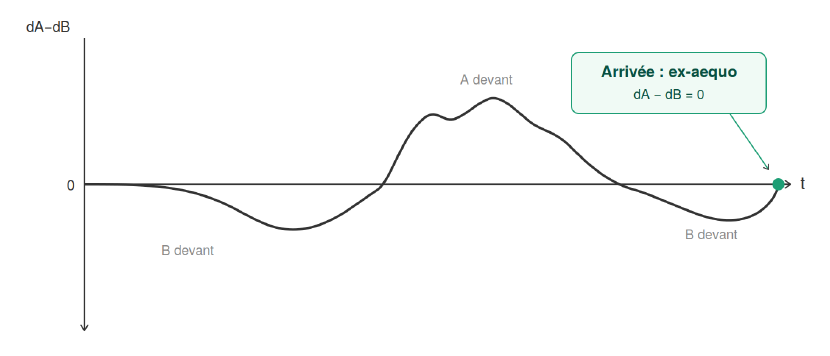}
\caption{The time at which both swimmers finish the race, producing a tie, according to Claude. The speech bubble reads ``Finish: tied, $d_A - d_B = 0$.'' \label{fig5}}
 \end{center}
\end{figure}

Claude's reply:  {\color{blue} ``You are absolutely right, I apologize! If the last point is exactly at $0$, then both swimmers finish at the same time---it is a tie. Here is the corrected graph:'' } And he produced the graph in Figure \ref{fig5},
and added {\color{blue} ``The correct conclusion is therefore: a tie! Both swimmers touch the wall at exactly the same moment.''}

I decided not to continue probing this incorrect answer for a few minutes and to come back to Question (a). I stated: ``You had noticed oscillations towards the middle of the graph, when A is in front. What do these oscillations mean?'' 

Claude's response:  {\color{blue} ``These oscillations correspond to the turns carried out in mid-race, when A is in front. Here is how they should be interpreted: 

 - Curve rises: $A$ is gaining on $B$. This corresponds to times where $A$ turns better than $B$---better tumble turn, pushes off faster, and increases his lead for a few moments.

  - Curve falls back: $B$ is catching up with $A$. This is the phase between two turns, where $B$, who is perhaps stronger in crawl, nibbles away the advantage that $A$ has just taken.
  
 In other words, these  oscillations reveal that $A$ is a better turner, but $B$ is a faster swimmer. If the race had been 50 m (no turns), $B$ would probably have dominated. It is the 200 m format, with its many turns, that  allows $A$ to stay competitive for so long.''}

 This sophisticated reply made me think of students who try to impress the examiner, without realizing that they are digging themselves deeper into a hole.

My next question: ``Can you show me on the graph the moment when $A$ begins his first turn?'' 

\begin{figure}[h]
\begin{center}
\includegraphics[height=4cm]{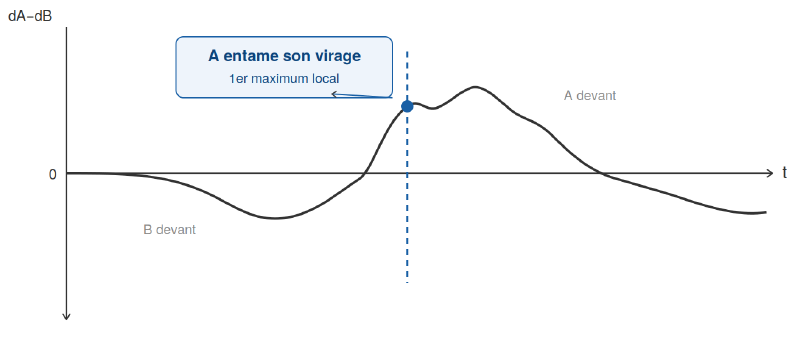}
\caption{The time where $A$ begins his first turn, according to Claude. The speech bubble reads ``$A$ starts his turn, $1^{\scriptsize \rm st}$ local maximum.'' \label{fig6}}
 \end{center}
\end{figure}

Claude replied  {\color{blue} ``The time when $A$ begins his turn 
corresponds to the first local maximum on the curve after the curve goes positive, just before the the first downwards portion in the oscillations''} and Claude produced the picture in Figure \ref{fig6} 

This was actually the first correct answer since the beginning of the conversation, so I tried to be helpful: ``And what is the moment when B begins his first turn?'' 

\begin{figure}[h]
\begin{center}
\includegraphics[height=4cm]{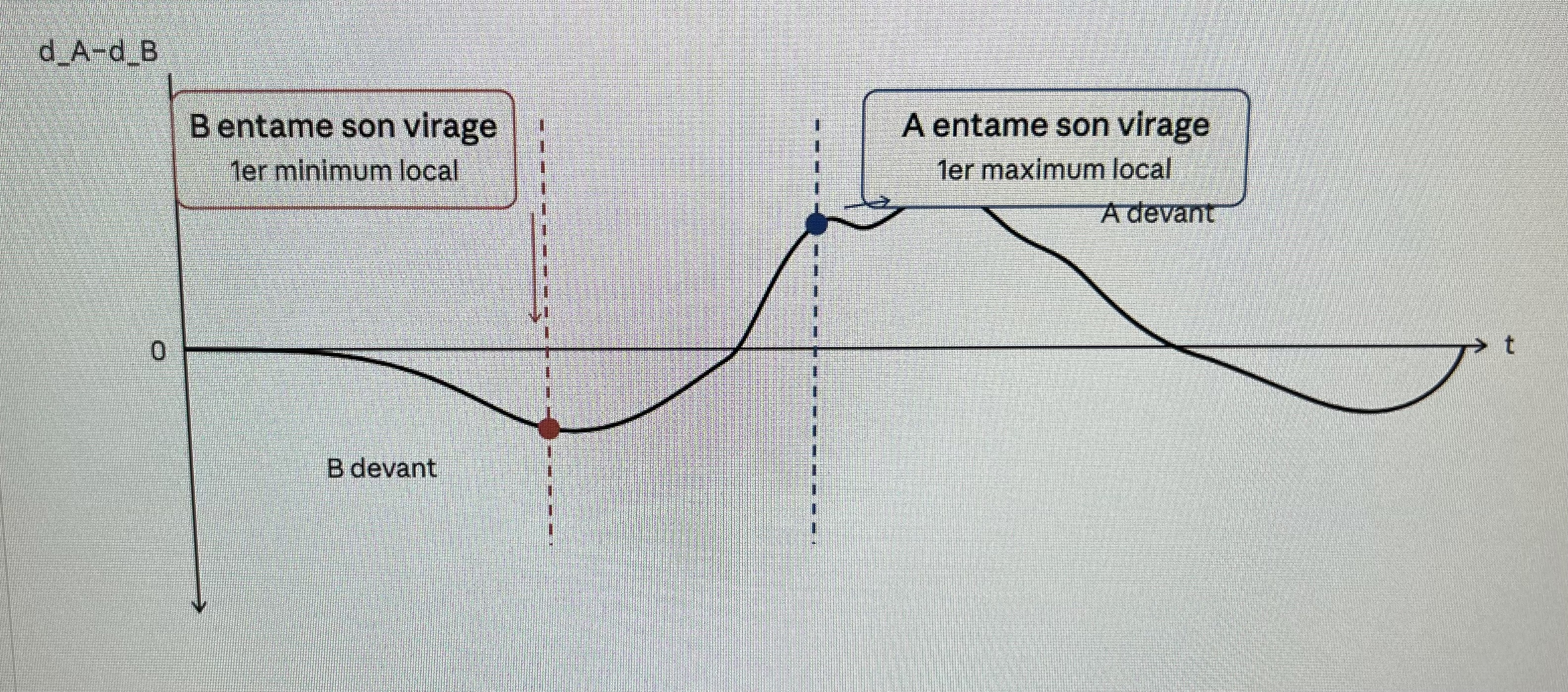}
\caption{Picture produced by Claude that purports to show where $B$ turned. The left speech bubble reads ``$B$ starts his turn, $1^{\scriptsize \rm st}$ local minimum,'' and the right bubble reads ``$A$ starts his turn,  $1^{\scriptsize \rm st}$ local maximum.'' \label{fig6a}}
 \end{center}
\end{figure}

Claude:  {\color{blue} ``Symmetric reasoning: Swimmer $B$ is in the lead at the beginning, therefore his first turn corresponds to the first local minimum of the graph---the moment when $B$ stops improving his lead on $A$ and begins his tumble turn. This is the point that we hade already identified together earlier on,''} and Claude produced the picture in Figure \ref{fig6a}\footnote[2]{For some unknown reason, Claude did not manage to give me the .pdf file with this figure. So I took a picture of the screen, and this is why the picture in Figure \ref{fig6a} looks different from the others. See \cite{claude}.}.

And Claude added {\color{blue} ``The general rule is quite elegant: turns of the leading swimmer correspond to local extrema—minimum when $B$ turns, maximum when $A$ turns.''}

 Again, a completely incorrect answer. So I asked ``Should A not turn very shortly after B?'' 
 
\begin{figure}[h]
\begin{center}
\includegraphics[height=4cm]{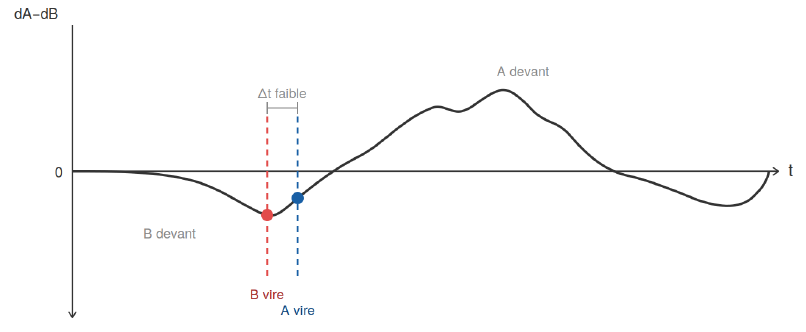}
\caption{The time where  $B$, then $A$, begin their first turn, according to Claude. \label{fig7}}
 \end{center}
\end{figure}

 Claude:  {\color{blue} ``Absolutely. Both swim the same distance, so $A$ necessarily turns shortly after $B$.
 On my earlier graph the two points were placed too far apart. 
 The first minimum and the first maximum are in fact very close in time, separated by a few seconds only: $B$ touches the wall, then $A$ touches it moments later, as in Figure \ref{fig7}.''}

  The fact that the tumble turns by $A$ and $B$ should occur close together in time sounded good to me, but the positions for these turns, as indicated by Claude in Figure \ref{fig7}, was completely off.

I decided to come back to the length of the race, which is in fact one of the most interesting questions and based on comparatively subtle features of the graph. So I told Claude: ``In fact, the graph is taken from a 100 m race in a 50 m pool. Would you like to modify some of your answers?'' 

Claude: {\color{blue} ``Yes, several things need to be revised:

    I had concluded that the race was 200 m by counting the oscillations, but it is 100 m—two lengths of 50 m,  therefore one turn per swimmer.
    
    With only one turn, the bumps cannot be turns. They must be natural pace variations during swimming---arm strokes, slight accelerations---or due to noisy measurements.


  Here is what remains valid:
    
   - $B$ is in the lead after the first length,

    - $A$ moves ahead after the turn,

    - they finish in a tie,

    - the first minimum is $B$'s tumble turn, the first maximum just afterwards is $A$'s tumble turn.''}

Since Claude's description of the oscillations was not was I was looking for, I asked ``In that case, what do the oscillations correspond to, given that there is no measurement noise?'' 

Claude: {\color{blue} ``Crawl is not a constant-speed motion—each arm stroke generates an impulse followed by a slight deceleration. With slightly different stroke frequencies or offset phases, $d_A - d_B$ oscillates with small periodic fluctuations around a mean trend. The oscillations reflect the arm stroke mechanics, not turns.''}

 Now anyone who has watched a swimming competition, particularly in crawl, knows that the swimmers progress very smoothly and there are no variations in speed that the eye can catch, except during turns. So I tried to give him a hint: ``According to your graph (Figure \ref{fig7}), the first turn appears at around 1/3 of the race. Is this normal?'' 
 
 Claude: {\color{blue} ``No, not at all. The turn must occur at roughly $T/2$ (after 50 m), where $T$ is the total race time. The second length is slightly slower (due to fatigue), so the turn is just below $T/2$, but never as early as $T/3$. My earlier graphs were poorly proportioned.

Is the minimum in your original graph placed around the midpoint of the time axis? That would be consistent with a properly represented 100 m race.'' }

\begin{figure}[h]
\begin{center}
\includegraphics[height=4cm]{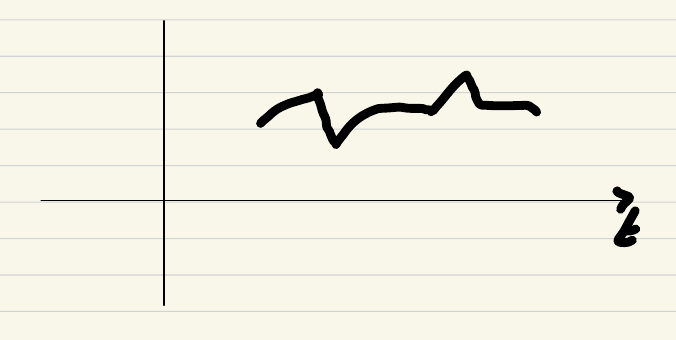}
\caption{A blow-up of the central part of the graph. \label{fig8}}
 \end{center}
\end{figure}

I consider that it is the examiner who should be asking the questions, not the entity being examined, so I did not reply, but I wanted to make sure that Claude had a clear picture of the oscillations that we were discussing. So I said ``Here is an enlargement of the central part of the curve,'' and I showed him Figure \ref{fig8}. 

Claude: {\color{blue} ``I see a very clear structure: a large deep minimum flanked by smaller oscillations. I reinterpret as follows: the deep minimum in the center is in fact the tumble turn---$B$ touches the wall, then $A$ touches it immediately afterwards. This is the quick drop in the curve, followed by a rise that is just as steep. This pronounced V-shaped form is typical of a turn: both swimmers decelerate, then push against the wall and go off again. 

   The small oscillations on each side of this minimum do correspond to the turning arms  that I mentioned: the irregular propulsion in crawl style swimming that makes the gap between the swimmers fluctuate.''
}

Claude does not seem to have noticed that the deep minimum that it is describing corresponds to the tumble turn of the lead swimmer only, and makes no mention of the nearby local maximum, which corresponds to the tumble turn of the second swimmer. But I decide not to bring this up for the moment.
   

 And Claude continues {\color{blue} ``This confirms that the turn occurs near the midpoint, which is consistent with a 100 m race in a 50 m pool; my previous graphs were poorly proportioned, with the minimum placed too early; the distinction between turn (large V) and  oscillations due to the swimming (small fluctuations) is clearly visible on this enlargement.'' }
 
 When Claude says ``the minimum placed too early,'' I am not sure what it means, since the minimum is a feature of the graph that is given, and Claude cannot choose where to place it. Claude's task is only to interpret what the location of the minimum may mean. And personally, I do not see any of the ``small fluctuations'' that seem to be ``clearly visible'' to Claude!

My next question: ``Is it not A who touches the wall first?'' 

Claude: {\color{blue} 
``Ah yes, you are right! I got confused. At the moment of the turn, the curve is positive ($d_A > d_B$), so $A$ is in the lead and touches the wall first, then $B$ a few instants later.

Here is what happens in the neighborhood of the deep minimum:

   1.~The curve drops abruptly: $A$ has just touched the wall and is doing his tumble turn, not advancing. $B$ keeps swimming and catches up, possibly even passes in front.
   
   2.~The minimum corresponds to the moment where the gap is maximal in $B$'s favor---$A$ still turning while $B$ is still swimming. 
   
    3.~The curve rises: $A$ pushes off, $B$ turns. $A$ is already in crawl and regains the lead.''}

Clearly, Claude's story does not match the curve: $A$ touches the wall, stops and turns, while $B$ is catching up. This corresponds to the steep drop in the curve. Then $A$ pushes off, so he regains some, but not all, of the time lost while turning, which is why after $A$'s turn, the curve is slightly lower than just before he started the turn. Then, a few second later, the steep upwards slope indicates that $B$ has stopped, makes his tumble turn and then, as $B$ pushes off, he regains some of the lost time, but not all. Claude seems to have missed these details completely.
 
 So I come back to the question: ``And in the end, who wins the race?'' 
 
 Claude: {\color{blue} ``According to the original graph, the curve ends exactly at zero—so a tie, as we had concluded. But honestly: the curve is freehand, and it is difficult to say with certainty whether the last point is exactly zero, slightly positive ($A$ wins), or slightly negative ($B$ wins). The precision of a freehand drawing does not allow a definitive conclusion.''}

Claude is correct that the graph is freehand, but whether the curve ends exactly at $0$, or very slightly above or below due to imprecision in the drawing, is in fact {\em not at all relevant} for determining who won the race.  Since Claude was again missing one of the most important features of the graph, I decided to help it a bit: ``At the end of the race, both swimmers have their hands on the side of the pool, so the difference between the distances is exactly equal to zero.'' 

\begin{figure}[h]
\begin{center}
\includegraphics[height=4cm]{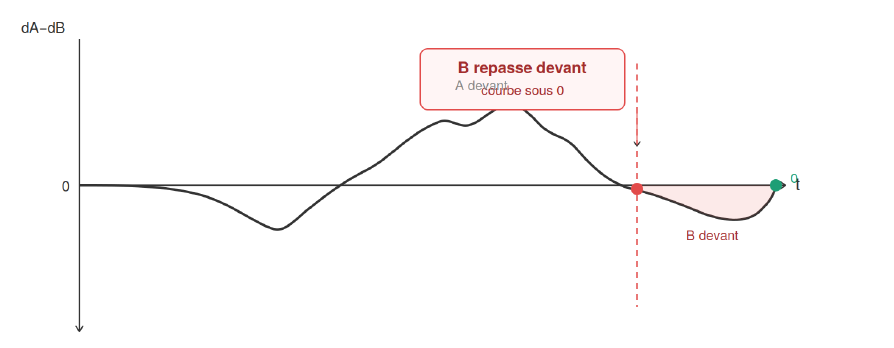}
\caption{Claude finally showing where $B$ takes the lead to win the race. The speech bubble reads ``$B$ passes in front again, the curve goes negative.'' Under the curve, it is written ``$B$ is ahead.'' \label{fig9}}
 \end{center}
\end{figure}

Claude: {\color{blue} ``Of course—I should have seen that immediately! Both cover exactly 100 m each, so $d_A(T) = d_B(T) = 100$ m, hence $d_A(T) - d_B(T) = 0$ by definition, no matter the precision of the drawing. 

The curve arrives at zero from below (negative just before the end), so $d_A < d_B$ in the final instants, meaning that $B$ reached 100 m first. $B$ wins the race.''} And to confirm its reasoning, Claude produced the picture in Figure \ref{fig9}.

And that is the end of the exam. With a substantial amount of prodding, Claude finally came to the correct conclusion.  Had Claude been a human being, this would have been a fairly painful exam, both for him and for me.

\section{Pierre's exam, then Catherine's}

Pierre's exam was more pleasant, probably for both of us. After reading Question (a), Pierre quickly noticed the bumps in the middle of the graph in Figure \ref{fig1} and realized that these must correspond to tumble turns. He had some trouble explaining why tumble turns would produce a graph of this kind, but he concluded correctly that this had to be a 100 m race.

Turning to Question B, Pierre had already understood where the turns were, and he pointed out that the curve was positive when the first swimmer begins his turn, and so he concluded correctly that it is Swimmer $A$ that is ahead at this time.

Concerning Question (c), Pierre initially thought that there was a tie between the two swimmers, because the curve ends at level $0$. So I pointed out, as I had for Claude, that at the end of the race, both swimmers have completed 100 m. With this hint, Pierre immediately realized (as had Claude when we finally got to this point), that since the curve is negative before it reaches $0$, $B$ has to be ahead. When I asked Pierre to show me the time on the graph at which $B$ wins the race, he had to think about it a bit, but then he realized that the beginning of the final segment that goes upwards must be when $B$ finishes the race and then $A$, who is close behind, also finishes shortly after $B$.
Finally, Pierre had no difficulty accurately recounting the race to a fictive radio audience.

Catherine's exam went similarly. Of course, since her primary training was not scientific, she did not know that the notation $d_A(t)$ shown on the vertical axis in Figure \ref{fig1} referred to the distance travelled by $A$ between times $0$ and $t$. But with light prodding, she did realize that the oscillations in the middle of the graph were related to the tumble turns, and therefore it was a 100 meter race. She noticed that towards the end of the race, the curve is negative, therefore $B$ is ahead, and she concluded that $B$ won. I did have to draw her attention to the fact that the curve necessarily ended up at $0$ before she could point to the moment when $B$ touches the side of the pool and wins the race. Finally, she concluded the exam by  recounting the race to the fictive radio audience in an entertaining way.

\section{Discussion and feedback to the ``examinees''} 

Let's continue imagining that Claude,  Pierre and Catherine are three students taking an exam with me, and that after the exam, I might give them some feedback on their performance.

For instance, if these conversations, with Claude, Pierre and Catherine, had occurred  during an oral exam, then it would have been clear to me that Pierre did much better than Catherine, and Catherine did better than Claude: Pierre was able to identify many of the key features of the graph (but not all) on his own. I would say that he ``connected the dots'' rather well, and this did require some concentration on his part: he even told me afterwards that he actually enjoyed this! Catherine had more difficulty, but still, she could ``put two and two together'' to understand what was going on. 

On the other hand, Claude was able to ``see all the dots,'' but he was mostly clueless about their meaning and ``unable to connect the dots.'' The graph in Figure \ref{fig3} clearly shows that he does not understand what he is talking about, even though Claude does state his answer with what I would interpret as confidence, were it a human being. 

With his Figure \ref{fig3}, Claude realizes that $B$ must be winning the race, but Claude changes the end of the graph so that it becomes meaningless. This last is not so important, because when drawing this graph, Claude may not be focussed on the end of the race, so I would not penalize it for this mistake. However, when later asked again who wins (Figure \ref{fig4}), Claude contradicts itself and concludes that the race is a tie, because it has noticed that the curve ends at $0$, but misinterprets this fact (as  Pierre and Catherine did  initially, too).

Claude initially states that the oscillations in the middle of the graph correspond to turns by the swimmers, which is correct, but when I ask where the first turn occurs, Claude points to a position about 1/3 of the way through the race (Figure \ref{fig3}), and when I express surprise at this fact, Claude again refers to the local minimum that it pointed out in Figure \ref{fig3}. Again, if Claude were human, then I would conclude that he was really clueless about what the graph is showing. 

Claude has a small confusion just after I show it Figure \ref{fig7}, asserting that $B$ touches the wall first, then switching its answer to $A$ when I suggest this. I would not penalize it for this switch: humans often say the opposite of what they really wanted to say, and then correct it when the mistake is pointed out.

Identifying the winner and realizing that the race was not a tie posed a problem to  Claude, Pierre and Catherine. And all three reacted correctly when I pointed out that at the end, the graph had to be at $0$ no matter who won the race.

I would not presume to compare the intelligence of Pierre with that of Catherine. I believe that both are very intelligent people, each in his/her own way. However, I am confident that in terms of {\em understanding} pre-calculus mathematics,  the ranking ``Pierre is better than Catherine and Catherine is better than Claude'' is accurately determined by this exam. Indeed, the comparison between Pierre and Catherine is already suggested by the simple fact that Pierre had a solid scientific training while Catherine studied literature. And the exam performances of these two individuals compared with that of the AI software speak for themselves.

Before attributing a grade to the three examinees, which is part of my job as professor, I should mention that at my university, we attribute grades between 1 and 6, with 6 being the best grade. The exam is passed if the grade is 4 or higher. I would give 5.5 to Pierre, 4.5 to Catherine, and 3 to Claude.

One could question the idea of attributing a grade based on a single sequence of questions on a fairly limited topic. However, for the last 30 years, when I have given an oral exam to a class of 30 students, say, I would prepare a list of 15 questions, each covering one or two topics that were covered during the semester-long course. Altogether, these 15 questions would cover essentially all the topics discussed in the course. An entering student would draw at random one question, which would then also be the question asked to the student entering the exam room immediately after. And the course grade would be determined by the student's answers to this single question. This method has its strengths and weaknesses, but I consider that during the 30 years that I have used it, this method has served my university well.

\section{Conclusions}

To be frank, I was amazed by Claude's ability to converse with me and to respond to my questions like an actual human being, even if this would be a human who is not very good at putting together basic pre-calculus concepts. There mere fact that Claude was able to ``see the graph,'' as given in a .pdf file, 
is in itself an amazing feat of programming. 
Claude identified all the important features of the curve, and I was very impressed by this: one could say that Claude {\em ``identified the location of all the dots,''} it just did not quite manage {\em ``to connect them''} and make use of the ``knowledge'' that it had extracted from the .pdf file. Whereas Pierre and Catherine, and particularly Catherine, who had only minimal training in mathematics, were able to use their innate intelligence to quickly understand the meaning of the graph. 

 I would say that Claude is almost capable of passing the famous Turing Test, which is a benchmark for machine intelligence proposed by the British mathematician and computer scientist Alan Turing in 1950: if a human judge cannot reliably tell if he is conversing with a man or with a machine, then the machine passes the test. I consider that Claude mostly passed this test, and that this extraordinary ability of Claude is a testimonial to the engineering know-how and ingenuity exhibited by all those who have contributed to make Claude and similar AI softwares what they are today. Ten years ago, or even five years ago, I would not have imagined that I would see the day when a computer software was capable of such a performance. As we say in French, to all these engineers and entrepreneurs: ``Chapeau bas~!''

For almost two years now, I have subscribed to a paying version of one of the main AI softwares available, most recently Claude, and I consider my \$20.-- per month payment to be a very good investment. AI does lots of things for me: it translates math documents from one language to another, it solves the math problems that I give to my students more quickly than I, and it answers my questions on essentially any topic (including questions that I asked it while preparing this article). Of course, one has to be careful: if the document that I ask an AI software  to translate is too long, then it sometimes translates the beginning and the end correctly, but puts unrelated words in the middle, possibly hoping that I will not notice, and maybe because it figures that for \$20 bucks per month, it should not be doing more for me than that. In this regard, it behaves somewhat like a tired and underpaid employee. And the solutions that it provides to math problems, while usually correct, look like they have been written by an average student, while I want to give my students fine-tuned solutions that look like they have been written by the professor (nevertheless, I feel that I save time by starting from Claude's solutions rather than from scratch, even if I still have to spend lots of time rewriting and optimizing Claude's solutions). And Claude gives helpful replies to most questions on most topics, but one really has to check the accuracy of the reply for any question where the answer is sensitive. 

Claude seems to have excellent general knowledge, but since all that it ``knows'' is presumably drawn from the Internet, there are still significant limitations. For example, after working for 40 years as a professional mathematician, there are many things that I know but have never written down in books or other publications of mine. Similarly, an experienced doctor will know much more about his specialty than he/she has published, so if I need a specialized medical opinion, I will definitively prefer to see a specialist than to converse with an AI, even if the AI has read every medical paper that was ever published. Of course, the day that computers will be able to read directly the information that is in our brains (and this day may not be so far away), the situation might change!

\begin{figure}[h]
\begin{center}
\includegraphics[height=4cm]{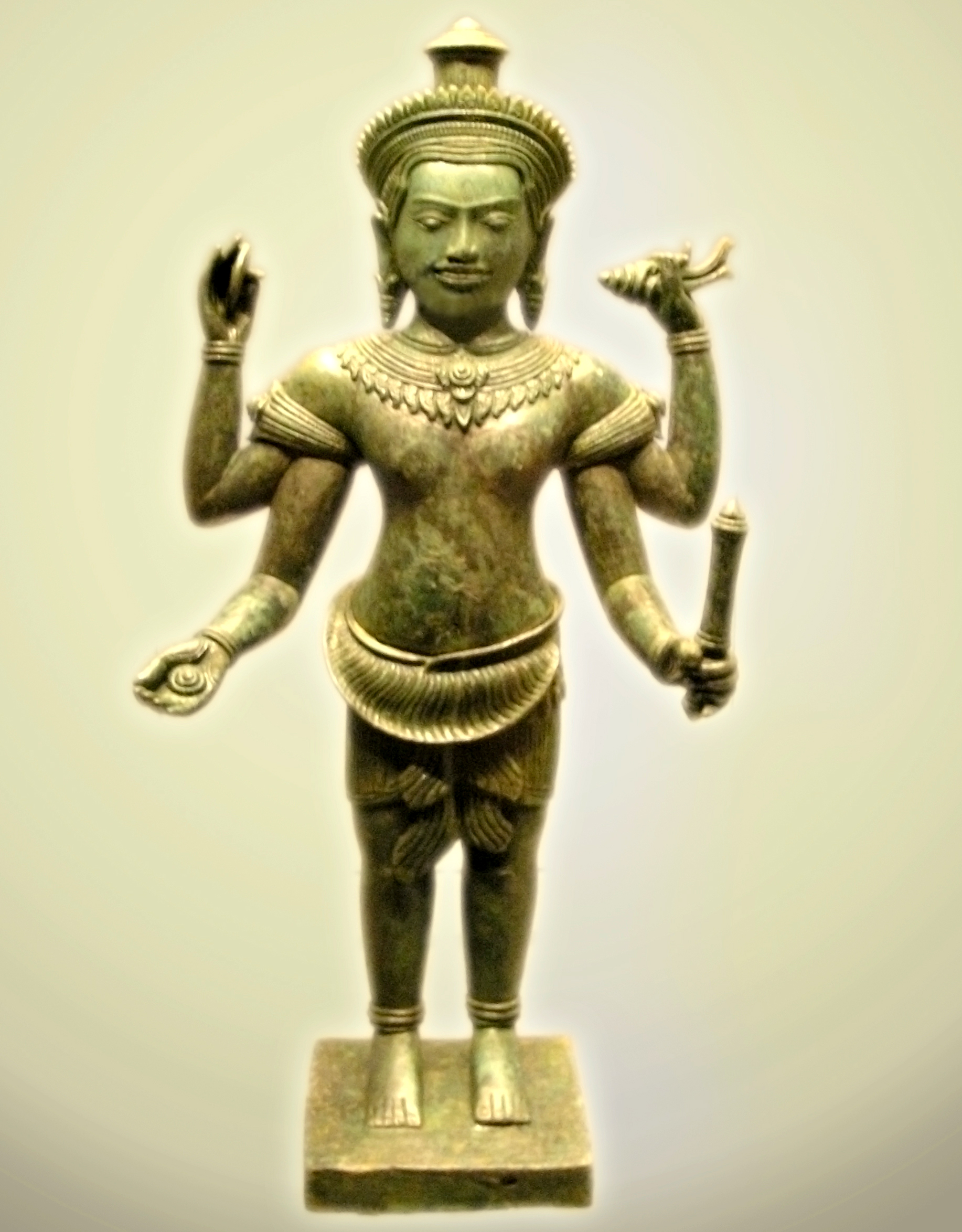}
\caption{The Hindu god Vishnu. Credit: Gryffindor, Public domain, via Wikimedia Commons. \label{fig10}}
 \end{center}
\end{figure}

I read in the news stories about how AI is helping research mathematicians prove new theorems. I have the utmost respect for the mathematicians involved in these activities, and I can well imagine how AI can help with mathematical research, because access to AI software is somewhat comparable not only to having at your fingertips the knowledge contained in all the books in the world, but also having at your service an experienced coder of computer software, an excellent debugger of LaTeX files, etc., in other words, an employee with numerous competences somewhat like the many-armed Hindu god Vishnu (Figure \ref{fig10}).
But it is not clear to me that the stories reported in the media, and particularly, in social media, accurately reflect the interactions of these distinguished mathematicians with AI software. And based on my own experiences with AI and on the fact that Claude failed miserably at a simple high-school math question, I do not see AI making original mathematical contributions in the near future.

AI is progressing very quickly, and it has to. Indeed, had it been possible to train an  AI software like Claude before smartphones were invented, then this software would be obsolete today, unless there had been a continuous effort to update it with new information. And the AI software of today will also have to be continuously updated, presumably, for the time being, with the new ideas and pieces of information generated by humans, if it is still to be relevant 10 or 20 years from now. For instance, having recorded and analyzed all the weather patterns everywhere in the world over the last 50 years is quite useful to make weather forecasts: an AI prediction then only requires looking for similar patterns in the past and using them to make predictions today. This is all well and good, but for instance, if the Gulf Stream one day changes direction, which is a possibility, given current global warming trends and such, then AI weather forecasts based on past weather data will be rendered useless. Whereas predictions based on physics and mathematical models will remain as useful as ever. So AI {\em is obliged} to progress, just so as not to become obsolete.

Current AI softwares are based on comparatively elementary mathematical algorithms, such as ``gradient descent,'' that we teach to undergraduate students. Generative AI uses slightly more sophisticated mathematical tools, based on properties of stochastic processes that are discussed in Master level courses, in order to produce images and videos from textual descriptions. As I stated above, one can only admire the engineering know-how and ingenuity that it takes to produce, from this theoretical framework, software as complex as Claude. Nevertheless, the fundamentals of life have an incredibly higher level of complexity than any existing computer code. For instance, the manner in which DNA polymerase reads the genetic sequence during replication is governed by purely random motion provoked by the impacts of fluid molecules on the DNA molecule, and the same type of random impacts, as studied by Albert Einstein \cite{einstein}, helps proteins fold and unfold and enzymes to meet their substrates in the cell. And it is not far-fetched to imagine that creativity, conscience, feelings and ideas are due in part to effects from quantum physics that occur in the human brain at the atomic or even sub-atomic level. So I do not believe that AI will exhibit any creativity on the level of Mozart, Einstein or Picasso for the foreseeable future.

This being said, having observed the path taken by chess and Kasparov, and by the game of Go and the AI software AlphaGo, I hope that Pierre and Catherine will not end up being the last human beings who, when pitted against one of the leading AI softwares of their day, were able to score a win for all humanity!

\noindent{\bf Author information:}
\medskip

\noindent{\scshape
Robert C.~Dalang,} Institut de Mathématiques,
\'Ecole Polytechnique Fédérale de Lausanne (EPFL),
CH-1015 Lausanne,
Switzerland.
Email: robert.dalang@epfl.ch

{\small
\begin{thebibliography}{AA}

\bibitem{claude} Dalang, R.C. A conversation with Claude on May 17, 2026. 10pp. Available on the website https://drive.switch.ch/index.php/s/KbD3kKcxvDwX7Rr

\bibitem{einstein} Einstein, A. Über die von der molekularkinetischen Theorie der Wärme geforderte Bewegung von in ruhenden Flüssigkeiten suspendierten Teilchen, {\em Annalen der Physik,} vol.~{\bf 17} (1905), pp. 549–560

\end{thebibliography}
}
\end{document}